\documentclass[letterpaper, 10 pt, conference]{ieeeconf}  

\IEEEoverridecommandlockouts                              
\usepackage{graphicx} 
\usepackage{amsmath} 
\usepackage{amssymb}  
\usepackage{cleveref}

\DeclareMathOperator*{\argmin}{argmin}
\newcommand{\eps}{\varepsilon}
\newcommand{\R}{\mathbb R}
\newcommand{\abs}[1]{\left \lvert #1 \right \rvert}
\renewcommand{\o}{\omega}

\title{\LARGE \bf
A Hamilton-Jacobi Formulation for Time-Optimal Paths of \\ Rectangular Nonholonomic Vehicles
}

\author{Christian Parkinson, Andrea L. Bertozzi, Stanley Osher
\thanks{Christian Parkinson is a graduate research assistant in the Department of Mathematics, UCLA,
        Los Angeles, CA, 90095 USA
        {\tt\small chparkin@math.ucla.edu}}%
	\thanks{Andrea L. Bertozzi and Stanley Osher are with the Faculty of the Department of Mathematics, UCLA {\tt\small bertozzi@math.ucla.edu, sjo@math.ucla.edu}}
}%

\begin{document}

\maketitle
\thispagestyle{empty}
\pagestyle{empty}

\begin{abstract}

We address the problem of optimal path planning for a simple nonholonomic vehicle in the presence of obstacles. Most current approaches are either split hierarchically into global path planning and local collision avoidance, or neglect some of the ambient geometry by assuming the car is a point mass.  We present a Hamilton-Jacobi formulation of the problem that resolves time-optimal paths and considers the geometry of the vehicle. 

\end{abstract}

\section{INTRODUCTION}

As autonomous vehicle technology becomes more and more prevalent, it is important to develop robust and widely applicable trajectory planning algorithms. Many such vehicles---planetary exploration rovers \cite{Tompkins}, flying drones \cite{Nieto}, or remote-controlled submarines \cite{Aguiar}---are subject to motion constraints which are \emph{nonholonomic}, depending not only on the configuration, but the velocity of the vehicle. Accordingly, much effort has been devoted to trajectory planning for general nonholonomic mechanical systems \cite{Colombo,Grushkovaskaya,Varricchio}.

\begin{figure}[b!]
\centering 
\includegraphics[width = 0.3\textwidth]{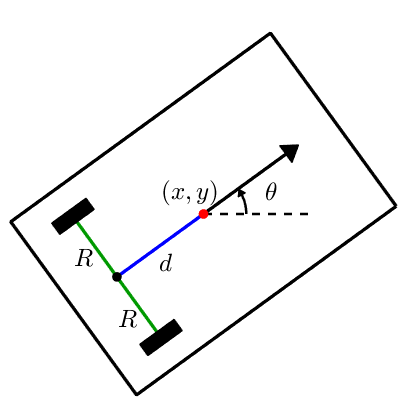}
\caption{A simple self-driving car.}
\label{fig:robotPic}
\end{figure}

One important example of a nonholonomic vehicle is a simple self-driving car. To track the motion of such a car, we model the current configuration using variables $(x,y,\theta)$: the spatial coordinate $(x,y)$ is the position of the center of mass of the car, and the orientation $\theta$ is the angle between the rear wheels and the horizontal, increasing in the counterclockwise direction. The car drives using actuators attached to the rear wheels that supply torque to each wheel individually, and steers using some mechanism separate from the rear wheels. The car has a rear axel of length $2R$, and a distance $d$ between the center of the rear axel and the center of mass, as pictured in figure~(\ref{fig:robotPic}). The motion of the car is constrained by a minimum turning radius, or equivalently a maximum angular velocity. This bound could be resolved in terms of $d,R$ and other parameters inherent to the steering mechanism. 

\subsection{Previous Work}
The problem of path planning for simple self-driving cars goes back to Dubins \cite{Dubins} who envisioned a vehicle that could move forward along paths constrained by a minimum turning radius. Later, Reeds and Shepp \cite{ReedsShepp} generalized the Dubins car to one that could also reverse direction. In both these cases, the problem was analyzed in a geometric and combinatorial fashion, discretizing the path into regions of straight-line movement and arcs of circles. The paths were designed to minimize length, and no obstacles were considered. Barraquand and Latombe \cite{Barraquand} added obstacles to the model, and devised a method of growing a reachability tree outward from the desired final configuration. Based on similar analysis, Agarwal and Wang \cite{AgarwalWang} assumed polygonal obstacles and presented an efficient algorithm for resolving paths that are robust to perturbation and nearly optimal. 

Since then, there has been increased effort to resolve optimal trajectories for such cars (and similar robots) using methods rooted in optimization and control.    In this case, the nonholonomic constraint is \begin{equation}\label{eq:nonholonomic}
\dot y \cos \theta - \dot x \sin \theta = d\dot \theta 
\end{equation} which ensures rolling without slipping and motion in the direction parallel to the rear wheels \cite{Barraquand,Fierro}. There have been several discrete and variational models of motion planning for these vehicles \cite{Fierro,Ali, Deluca, Shukla, Wu}. One advantage of models based on optimization is that they can seemlessly account for paths that are not only time-optimal, but consider energy consumption as well \cite{Duleba, KhanEtAl,WangEtAl,Verscheure}. Discrete models of this sort are often hierarchical, relying on a global path planner and a local collision avoidance algorithm \cite{AlonsoMora, Lee2017}.

A model for curvature constrained motion based on dynamic programming and a Hamilton-Jacobi-Bellman equation was introduced by \cite{TakeiTsai1,TakeiTsai2}, where obstacles are included, but the car is simplified to a point mass, meaning extra concern is required near obstacle boundaries. Similar Hamilton-Jacobi type models for optimal path planning in other contexts are quite common, and include level set methods \cite{Parkinson} and fast-marching methods \cite{Tsitsiklis,Sethian1}.

\subsection{Our Contribution}

We present a model for optimal path planning of nonholonomic self-driving cars based on a Hamilton-Jacobi formulation, but considering the geometry of the vehicle. Our approach is akin to that of \cite{TakeiTsai1,TakeiTsai2}. However, those authors simplify the car to a point mass and accordingly, must either create a buffer region around an obstacles \cite{TakeiTsai1} or opt for a semi-Lagrangian path planning approach \cite{TakeiTsai2}. If we do not make the simplification, we can maintain the Hamilton-Jacobi approach. The Hamilton-Jacobi formulation has the natural advantage that it averts the need for hierarchical planning algorithms. Additionally, this approach can provide optimal trajectories from all starting configurations to a desired final configuration, as opposed to variational methods which typically resolve a single trajectory. The general steady-state Hamilton-Jacobi equation takes the form \begin{equation} \label{eq:HJ} H(x,\nabla u) = f(x).\end{equation} Because the equation is nonlinear, special care is needed to solve Hamilton-Jacobi equations numerically. Accordingly, we present an upwind sweeping scheme that traces the characteristics outward from a desired final configuration. 

\section{The Hamilton-Jacobi Formulation}

The Hamilton-Jacobi formulation for optimal-path planning is based on the dynamic programming principle \cite{Bellman}. One begins with a controlled equation of motion, and derives a partial differential equation satisfied by the value function. For our motion, we consider a kinematic equation that neglects some of the dynamics, but is sufficient for our purposes \cite{Wu}.

\subsection{Kinematics \& Control Problem}

Assume the car moves throughout a domain $\Omega \subset \R^2$ that is segmented into free space and obstacles: $\Omega = \Omega_{\text{free}} \cup \Omega_{\text{obs}}$. The current configuration of the car is given by $(x,y,\theta) \in \Omega \times [0,2\pi)$ as described above. As the car moves, it obeys the equations \begin{equation}\label{eq:motion} \begin{split} \dot x &= v \cos \theta - \o W d \sin \theta, \\ \dot y &= v \sin \theta + \o W d \cos \theta, \\ \dot \theta &= \o W, \end{split} \end{equation}  where $W>0$ is the maximum angular velocity, which bounds the curvature of a path, and $(v, \o) \in [-1,1]$ are the control variables representing tangential and angular velocity, respectively. This model assumes instantaneous changes in $(v,\o)$ which is akin to assuming infinite acceleration; this is what we mean when we say we are neglecting some of the actual dynamics. For any configuration $(x,y,\theta) \in \Omega \times [0,2\pi)$ let $C(x,y,\theta) \subset \R^2$ denote the space occupied by the car. The shape could be arbitrary, but for our car, this will be a rectangle of height $2R$ and width $2d$ that is centered at $(x,y)$ and then rotated by $\theta$. We call a configuration \emph{admissable} if $C(x,y,\theta) \cap \Omega_{\text{obs}} = \varnothing$. Next, suppose we are given a desired final configuration $(x_f,y_f,\theta_f)$ which is admissable. We call a trajectory $(x(\cdot),y(\cdot),\theta(\cdot))$---defined for $t\in[0,T]$---\emph{admissable} if it obeys \eqref{eq:motion} for all $t \in (0,T]$, $(x(t),y(t),\theta(t))$ is an admissable configuration for all $t \in [0,T]$, and $(x(T),y(T),\theta(T)) = (x_f,y_f,\theta_f)$. Given a starting point $(x,y,\theta)$, the goal is to choose $(v(\cdot),\o(\cdot))$ so as to minimize travel time $T$ among all admissable trajectories.

\subsection{Value Function \& Hamilton-Jacobi-Bellman Equation}

Denote by $\mathcal A(x,y,\theta,T)$ the set of admissable trajectories beginning at the configuration $(x,y,\theta)$ and requiring time less than $T$ to traverse. We define the travel-time function: \begin{equation}\label{eq:valueFunc} u(x,y,\theta) = \inf \{T \, : \, \mathcal A (x,y,\theta,T) \neq \varnothing \},\end{equation} where the infimum is taken over the control values $(v(\cdot),\o(\cdot))$. We formally derive a partial differential equation solved by $u$. The dynamic programming principle tells us that if $\delta > 0$, then \begin{equation}\label{eq:DPP} u(x(t),y(t),\theta(t)) = \delta +\inf\{u(x(t+\delta),y(t+\delta),\theta(t+\delta))\},\end{equation} where now the infimum is taken over the values $(v(\cdot),\o(\cdot))$ on the interval $(t,t+\delta)$. Intuitively, equation \eqref{eq:DPP} states that traveling optimally for time $\delta$ will decrease the remaining travel time by exactly $\delta$. This is an expression of the fact that globally optimal paths are also locally optimal. Assuming $u$ is smooth, we can divide by $\delta$ and send $\delta \to 0$ to see \begin{equation}\label{eq:intermediate}
-1 = \inf\left\{\dot x u_x + \dot y u_y + \dot \theta u_\theta\right\},
\end{equation} whence \eqref{eq:motion} yields \begin{equation}\label{eq:HJB1}
-1 = \inf_{v,\o} \left\{\begin{split} & (u_x\cos \theta + u_y \sin \theta)v \,\,+ \\& \,\,\,\,\,\,\,W(-du_x\sin\theta + du_y \cos\theta + u_\theta)\o \end{split}\right\}. \end{equation} The infimum can be resolved explicitly, showing that the optimal control values are given by  \begin{equation} \label{eq:controls} \begin{split} v &= -\text{sign}(u_x \cos \theta + u_y \sin \theta), \\
\o &= -\text{sign}(-du_x\sin \theta + du_y \cos \theta + u_\theta), \end{split} \end{equation} where $u(x,y,\theta)$ solves the Hamilton-Jacobi-Bellman equation 
\begin{equation} \label{eq:HJB} 
\left.\begin{split} 1=& \,\lvert u_x\cos\theta + u_y \sin\theta\rvert + \\& \,\,\,\,\,\,\,\,\,\,\,W \lvert-du_x\sin\theta + du_y\cos \theta + u_\theta\rvert. \end{split}\right.
\end{equation} The above computation only holds rigorously when $u$ is smooth which may not be the case. However, as long as $u$ remains continuous it will be the unique viscosity solution of \eqref{eq:HJB} \cite{CrandallLions,Bardi1997}, and if $u$ is discontinuous, one can still maintain existence and uniqueness by passing to a yet weaker notion of solution \cite{TakeiTsai2}. The travel time from the final configuration to itself is zero, so we impose the ``boundary" condition $u(x_f,y_f,\theta_f) = 0$. Similarly, to account for obstacles we set $u(x,y,\theta) = +\infty$ for any $(x,y,\theta)$ that is not admissable. 

Note \eqref{eq:controls} results in a bang-bang controller: $v,\o = \pm 1$. Along trajectories where the car is already oriented in the correct direction, we find $-du_x\sin\theta + du_y\cos \theta + u_\theta = 0$ so that $\o = 0$ is also a possibility. This agrees with the analysis in \cite{Dubins,ReedsShepp,Barraquand} where it is proven that optimal paths consist of straight lines and arcs of circles of minimum radius.

\section{Numerical Methods}

Since Hamilton-Jacobi equations are nonlinear and have solutions that develop kinks, some care is needed when solving them numerically. We would like to develop a sweeping scheme similar to those in \cite{TakeiTsai2,TsaiOsherSweep}. The primary concerns for such a scheme are that it should be upwind and monotone. For a general discussion of numerical analysis of Hamilton-Jacobi equations see  \cite{CrandallLions, CrandallMajda, OsherShu, Oberman}.

\subsection{Upwind Sweeping Scheme for \eqref{eq:HJB}} 

For simplicity, we describe the numerics on a rectangular domain $\Omega = [a,b] \times [c,d]$. Fixing $I,J,K \in \mathbb N$, let $(x_i)_{i=0}^I, (y_j)_{j=0}^J, (\theta_k)_{k=0}^K$ be the uniform dicretization of $\Omega \times [0,2\pi]$ (so that, for example, $\Delta x = (b-a)/I$) and let $u_{ijk}$ be the numerical approximation to $u(x_i,y_j,\theta_k)$. To discretize \eqref{eq:HJB} in a fully upwind manner, define\begin{equation}
\label{eq:AB} \begin{split}
A_k(v,\o) &= v \cos \theta_k - \o Wd \sin \theta_k, \\
B_k(v,\o) &= v \sin \theta_k + \o Wd \cos \theta_k, \\
a_k(v,\o) &= \text{sign}(v \cos \theta_k - \o Wd \sin \theta_k), \\
b_k(v,\o) &= \text{sign}(v \sin \theta_k + \o Wd \cos \theta_k).
\end{split} 
 \end{equation} Rearranging \eqref{eq:HJB1} shows that \eqref{eq:HJB} is equivalent to  \begin{equation}\label{eq:HJBuxuyus}
-1 = \inf_{v,\o} \left\{ A_{k}(v,\o) u_x + B_k(v,\o)u_y + \o W u_\theta \right\}. \end{equation} The upwind approximation to each derivative term in \eqref{eq:HJBuxuyus} at $(x_i,y_j,\theta_k)$ is given by \begin{equation} \label{eq:upwindApprox} \begin{split}
\big( A_{k}(v,\o) u_x \big)_{ijk} &= \abs{A_k(v,\o)}  \left(\frac{u_{i+a_k(v,\o),j,k} - u_{ijk}}{\Delta x}\right), \\
\big( B_k(v,\o) u_y \big)_{ijk} &= \abs{B_k(v,\o)}  \left(\frac{u_{i,j+b_k(v,\o),k} - u_{ijk}}{\Delta y}\right), \\
( \o W u_\theta)_{ijk} &= \abs{\o}W  \left(\frac{u_{i,j,k + \text{sign}(\o)} - u_{ijk}}{\Delta \theta}\right).
 \end{split} \end{equation} For a particular pair $(v,\o)$, if we plug these approximations into \eqref{eq:HJBuxuyus}, we can solve for $u_{ijk}$ in terms of the values at neighboring nodes. This shows that \begin{equation} \label{eq:update} \begin{split} u^*_{ijk}(v,\o) =  \bigg(1 &+ \frac{A_k(v,\o)}{\Delta x} u_{i+a_k(u,v),j,k} \\ &+ \frac{B_k(v,\o)}{\Delta y} u_{i,j+b_k(v,\o),k} \\&+ \frac{\abs \o W}{\Delta \theta}u_{i,j,k+\text{sign}(\o)} \bigg) / \\ &\bigg(\frac{A_k(v,\o)}{\Delta x} + \frac{B_k(v,\o)}{\Delta y} + \frac{\abs \o W}{\Delta \theta}\bigg) \end{split} \end{equation} is an upwind, first-order approximation to the solution of \eqref{eq:HJB} when the pair $(v,\o)$ gives the correct control values at node $(i,j,k)$. Together with the boundary conditions $u(x_f,y_f,\theta_f) = 0$ and $u(x,y,\theta) = +\infty$ at inadmissable configurations,  this suggests a sweeping scheme of the form: \begin{enumerate} 
 \item[($i$)] \emph{Initialization.} Set $u^0_{ijk} = 0$ at the nodes closest to $(x_f,y_f,\theta_f)$, and $u^0_{ijk} = +\infty$ (or some large number) for all other nodes. For all $(x_i,y_j,\theta_k)$ that are not admissable, add $(i,j,k)$ to the set $\Omega^*_{\text{obs}}$.  
\item[($ii$)] For all $i = 1 : I-1$, $j = 1: J-1$, $k = 0:K-1$ with $(i,j,k) \not\in \Omega^*_{\text{obs}}$, compute $u^*_{ijk}(v,\o)$ according to \eqref{eq:update} using the values $u^{n-1}_{ijk}$.
\item[($iii$)] Set $u^n_{ijk} = \min(\min_{v,\o} u^*_{ijk}(v,\o),u^{n-1}_{ijk})$
\item[($iv$)] Repeat steps $(ii)$,$(iii)$, sweeping through the indices in alternating directions until all combinations of sweeping directions have been performed. (This should be $8$ total sweeps: [$i$-forward,$j$-forward,$k$-forward],[$i$-forward,$j$-forward,$k$-backward], etc.)
\item[($v$)] Repeat steps $(ii),(iii),(iv)$ until convergence. 
 \end{enumerate}
 
 \subsection{Implementation Notes} We include a few implementation notes regarding the sweeping scheme. First, in step $(ii)$, the minimum over $(v,\o)$ corresponds exactly to the infimum in \eqref{eq:HJB1},\eqref{eq:HJBuxuyus}. We then choose to update $u^n_{ijk}$ only if the new value is smaller than the previous value. This ensures that the scheme is monotone \cite{Oberman}. Second, there is no need for computational boundary conditions at $i=0,I$, $j = 0,J$. Those nodes are never updated, so their values will remain large; this ensures that the car stays in the computational domain, and the upwind nature of the scheme ensures that those values do not effect the solution at interior nodes (likewise, the values at inadmissable configurations are never updated but will not affect the values at nearby admissable configurations). At $k = K$, one should enforce a periodic boundary condition identifying these values with $k=0$. Third, sweeping is carried out in the Gauss-Seidel sense: updating nodes and then using the most recently updated values as you go.  Fifth, to slightly reduce computational load, one can pre-compute the values $A_k(v,\o), B_k(v,\o), a_k(v,\o),b_k(v,\o)$ since they are static during the iteration. Finally, one can test convergence in any number of ways. We suggest the criterion $\sup_{ijk} \lvert u^n_{ijk} - u^{n-1}_{ijk} \rvert < \eps$ for a specified tolerance $\eps$. 
 
\subsection{Resolving Optimal Trajectories}
 While it is not stated in the algorithm above, the control values can be established during the sweeping by setting \begin{equation} \label{eq:controlUpdate}(v^n_{ijk},\o^n_{ijk}) = \argmin_{v,\o} u^*_{ijk}(v,\o) \end{equation} when the minimum is used, and $(v^n_{ijk},\o^n_{ijk}) = (v^{n-1}_{ijk},\o^{n-1}_{ijk})$ when $u^n_{ijk} = u^{n-1}_{ijk}$. Alternatively, after computing the value function $u(x_i,y_j,\theta_k)$, one can interpolate to off grid values and compute $(v,w)$ using \eqref{eq:controls}. Under fairly mild conditions on the Hamiltonian $H$, one can prove that solutions to \eqref{eq:HJ} remain Lipschitz continuous; hence differentiable almost everywhere \cite{Bardi1997}. This means that the control values---and thus optimal trajectories from $(x,y,\theta)$ to $(x_f,y_f,\theta_f)$---are uniquely determined unless $(x,y,\theta)$ lies in a set of measure zero. Note, the only crucial piece of datum is the final configuration $(x_f,y_f,\theta_f)$. Once that is specified and we have solved \eqref{eq:HJB}, we can compute the optimal trajectory from any other configuration in the domain to the final configuration.

We can then generate the optimal trajectory from a given point $(x,y,\theta)$ to $(x_f,y_f,\theta_f)$ by setting $(x(0),y(0),\theta(0)) = (x,y,\theta)$ and integrating \eqref{eq:motion} until the time $T = u(x,y,\theta)$. For our purposes, we resolve $(v,\o)$ according to \eqref{eq:controls} by interpolating $u(x,y,\theta)$ to off grid values and using a centered difference approximation for $(u_x,u_y,u_\theta)$. We integrate \eqref{eq:motion} using forward Euler time-stepping.  

\section{Results \& Observations}

We test our algorithm using the spatial domain $\Omega = [-1,1] \times [-1,1]$. In our tests, we set $R = 0.04$, $d = 0.07$ and $W = 4$, meaning that the minimum radius of a circle that the car can traverse is $1/4$. Note these are all dimensionless parameters, used solely for testing purposes. In our simulations, we used a uniform discretization with 200 grid nodes in each direction. Depending on the simulation, the sweeping scheme required roughly 25 iterations, though it took longer in simulations with obstacles. This number could likely be decreased by introducing more accurate finite difference approximations as suggested by \cite{TakeiTsai1}.

\begin{figure}[b!]
\centering
	\includegraphics[width=0.45\textwidth]{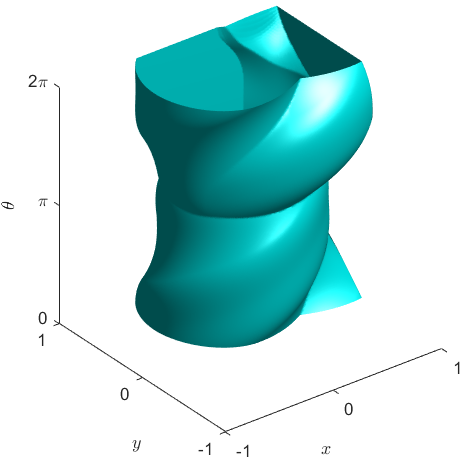}
	\caption{The isocontour  $u(x,y,\theta) = 1$, $(x_f, y_f, \theta_f) = (\tfrac 1 2, \tfrac 12, 0)$.}
	\label{fig:uLevelSets3d}
\end{figure}

\begin{figure}[b!]
\includegraphics[width=0.45\textwidth]{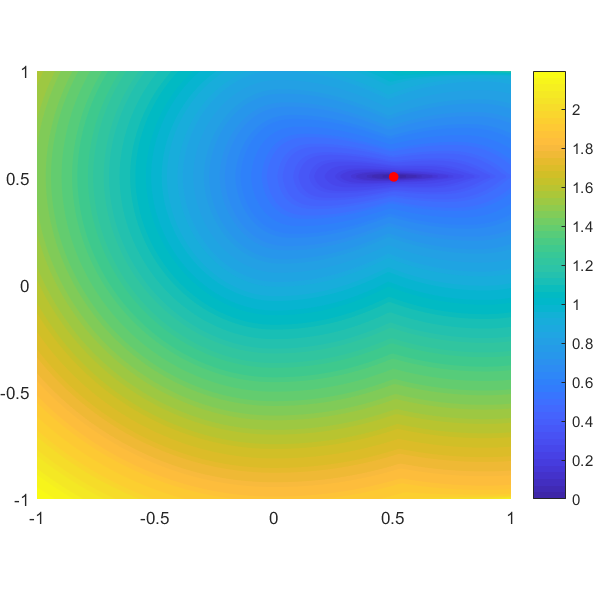}
\caption{The contour map of $u(x,y,0)$, $(x_f, y_f, \theta_f) = (\tfrac 1 2, \tfrac 12, 0).$}
\label{fig:uLevelSets}
\end{figure}

We first computed the value function $u(x,y,\theta)$ for the final orientation is $(x_f, y_f, \theta_f) = (\tfrac 1 2, \tfrac 12, 0)$ with no obstacles. This function $u(x,y,\theta)$ gives the optimal travel time from $(x,y,\theta)$, given that the car must end at the point $(x_f, y_f) = (\tfrac 1 2, \tfrac 1 2)$ facing horizontally in the positive $x$-direction. \Cref{fig:uLevelSets3d} displays the isocontour $\{(x,y,\theta) \in \Omega \times [0,2\pi] \,\, : \,\, u(x,y,\theta) = 1\}$. One interesting note here is the approximate symmetry across the line $\theta = \pi$, which is an expression of periodicity in the value function. When $d = 0$, we do indeed have $\pi$-periodicity: $u(x,y,\theta) = u(x,y,\theta+\pi)$ \cite{TakeiTsai2}. When $d \neq 0$, this is only approximate. \Cref{fig:uLevelSets} shows a contour map of the function $u(x,y,0)$, which gives the travel time from different starting positions if the car is already facing in the positive $x$-direction. The final position $(x_f,y_f) = (\tfrac 1 2, \tfrac 1 2)$ is represented by the red dot. As a sanity check, we note that along the line $y = \tfrac 1 2$, the value is given by $u(x,\tfrac 1 2,0) = \abs{x - x_f}$ since the optimal path merely includes pulling forward or reversing into the spot. 

\begin{figure}[t!]
\centering
	{\setlength{\fboxrule}{2pt}
\setlength{\fboxsep}{0pt}
\fbox{\includegraphics[width=0.48\textwidth,trim=80 50 50 50, clip]{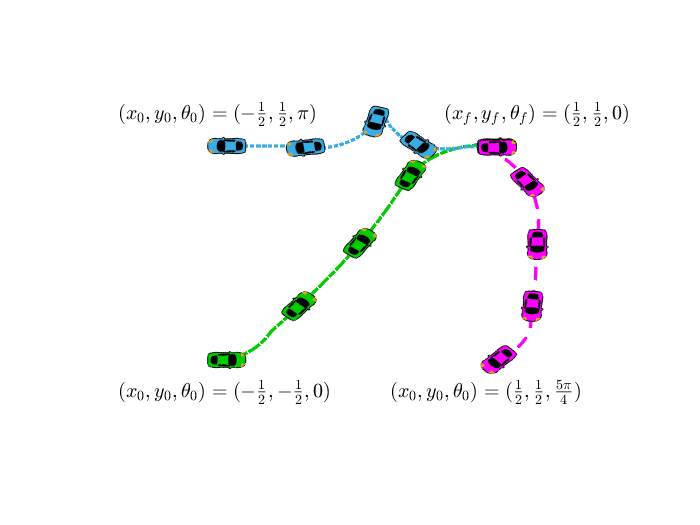}}}
	\caption{Three paths from different starting configurations to the final configuration $(x_f,y_f,\theta_f) = (\tfrac 1 2, \tfrac 1 2, 0)$.}
	\label{fig:threePaths}
\end{figure}

Using this value function, we can compute optimal trajectories from any initial configuration to the final configuration. \Cref{fig:threePaths} displays optimal paths originating from three different initial configurations. The initial and final configurations are labeled on the plot. The blue, green and pink car icons represent the vehicles, and the orientation is given by the direction the headlights are facing. \Cref{fig:threePaths} shows that the best strategies for the blue and pink car involve traveling large portions of the path in reverse ($v=-1$), before pivoting and achieving the final configuration while moving in the forward direction. By contrast, the green car only travels in the forward direction ($v= 1$). 

\begin{figure}[t!] 
\centering
	{\setlength{\fboxrule}{2pt}
\setlength{\fboxsep}{0pt}
\fbox{\includegraphics[width=0.48\textwidth,trim=80 50 50 50, clip]{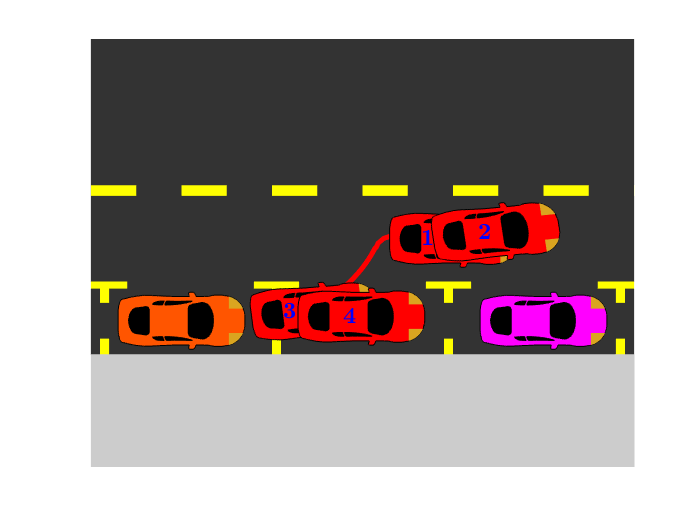}}}
\caption{A car parallel parking demonstrates an optimal path with two kinks. The numbers denote successive positions along the path.}
\label{fig:parallelPark}
\end{figure}

As stated above, Reeds and Shepp \cite{ReedsShepp} analyzed this problem in the case that $d = 0$ so that the car is a point mass. They proved that the optimal path between two points consists of a finite number of straight lines and arcs of circles of minimum radius. Further, they proved that while kinks will occur as the car switches driving direction, the optimal path requires no more than two kinks. Our simulations empirically confirm this; in the examples in \cref{fig:threePaths}, none of the paths required more than one kink. For an example of an optimal path with two kinks, consider the parallel parking problem displayed in \cref{fig:parallelPark}. In this example $(x_0,y_0,\theta_0) = (x_f+2d,y_f + 3R,\theta_f)$ and the car is plotted at four points along the path: the initial position, the two kinks, and the final position. 

\begin{figure}[b!]
\centering
	{\setlength{\fboxrule}{2pt}
\setlength{\fboxsep}{0pt}
\fbox{\includegraphics[width=0.48\textwidth,trim=80 50 50 50, clip]{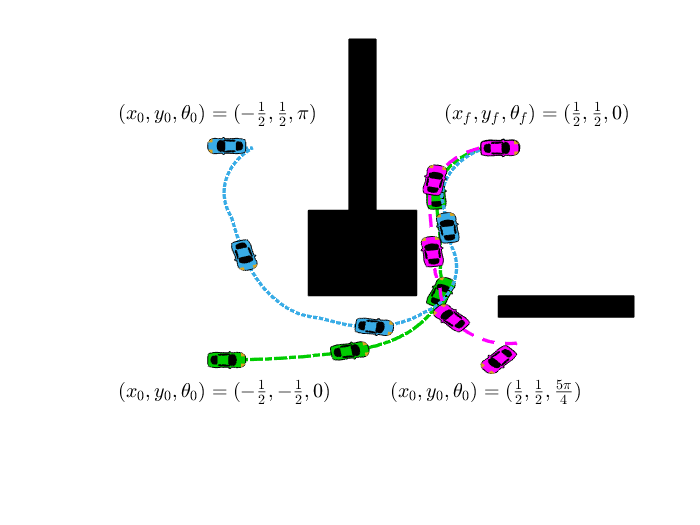}}}
	\caption{Three paths from different starting configurations to the final configuration $(\tfrac 1 2, \tfrac 1 2, 0)$, with obstacles [black].}
	\label{fig:threePathsObs}
\end{figure}

\begin{figure}[b!]
\centering
{\setlength{\fboxrule}{2pt}
\setlength{\fboxsep}{0pt}
\fbox{\includegraphics[width=0.48\textwidth,trim=100 62.5 62.5 60, clip]{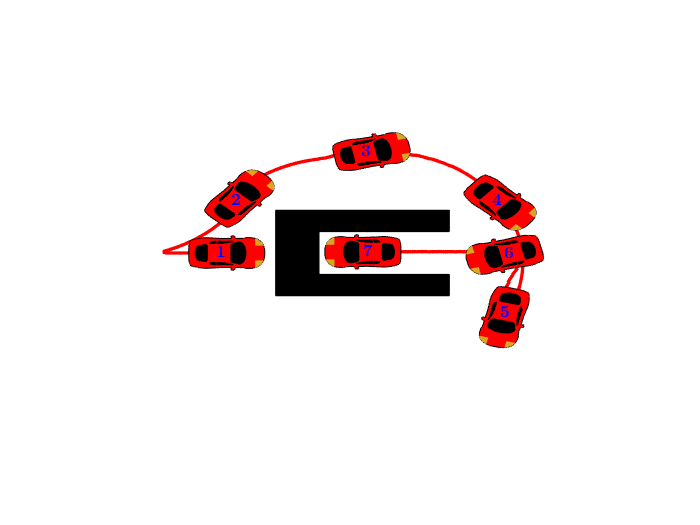}}}
\caption{A car parking in a very narrow spot. The numbers denote successive positions along the path.}
\label{fig:parkingSpot}
\end{figure}

Lastly, we introduce obstacles. To reiterate, the algorithm for solving \eqref{eq:HJB} is the same, except that the value function is not updated at nodes corresponding to illegal configurations---those which would cause the car to collide with an obstacle. In \cref{fig:threePathsObs}, we compute the optimal paths from the same three starting configurations as in \cref{fig:threePaths} but now with obstacles [black] hindering the cars' movement. In \cref{fig:parkingSpot}, we have a car pulling into a very narrow parking spot. Note that no extra consideration (in the form of local collision avoidance) was necessary to resolve this path. Here the width of the parking spot is only $0.1$ and the width of the car is $2R=0.08$. Thus if we buffered the obstacles, the final configuration would likely by illegal, and if we approximated the car by a point mass, it would likely take on some illegal configurations. As an aside, the optimal trajectory into the parking spot has three kinks, showing that the result of \cite{ReedsShepp}---stating that only two kinks are sufficient---is not true in the presence of obstacles. Optimal paths with polygonal obstacles are also considered by \cite{Barraquand}.

\section{Conclusion and Future Work}

In this manuscript, we presented a Hamilton-Jacobi-Bellman formulation for optimal path planning of nonholonomic vehicles, accounting for impassable obstacles and the actual geometry of the vehicle. We developed an upwind sweeping scheme to solve the Hamilton-Jacobi-Bellman equation for the value function. We validated our model in the presence and absence of obstacles and compared with the classical results for curvature constrained motion. We note that no extra considerations were needed when dealing with obstacles since the geometry of the car is not being neglected.

We propose two avenues for future work. First, rather than only considering time-optimal paths, one can very easily incorporate energy optimization into the this formulation. Indeed, one form of marginal energy cost along the path takes the form $E(v,\o) = \frac 1 2 m v^2 + \frac 1 2 I \o^2$ where $m$ is the mass and $I$ is the moment of inertia of the vehicle \cite{Fierro}. Thus we could change the cost functional to something like \begin{equation} \label{eq:energyCost} 
C[v(\cdot),\o(\cdot)] = \int^T_0 \big[\lambda + (1-\lambda)E(v(t),\o(t))\big]dt
\end{equation} for some weight $\lambda \in [0,1]$. In the $\lambda = 1$ case, the cost is merely the travel time and we revert to the model presented here. Accounting for energy does not complicate the model any further, but the numerics become more difficult since it is no longer a bang-bang control problem, and it is likely impossible to devise a simple update for a sweeping scheme. 

Second, curvature constrained motion can be easily extended to higher dimensions. In three spatial dimensions, one can consider different types of curvature constraints to account for vehicles like airplanes or submarines. These cases are also amenable to a Hamilton-Jacobi formulation, though they will be higher dimesional since they require dynamic models and hence second order controllers. Thus grid-based numerical methods will run into the curse of dimensionality and one may need to employ non-grid based methods like those presented in \cite{Lin}. Other modeling concerns would need to be addressed as well. For example, airplanes must maintain a minimum cruising velocity, and the braking mechanism for an airplane or submarine will be vastly different than that of a car. 

\section{ACKNOWLEDGMENT}

This research is funded by an academic grant from the National Geospatial-Intelligence Agency (Award No. \#HM0210-14-1-0003, Project Title: Sparsity models for spatiotemporal analysis and modeling of human activity and social networks in a geographic context). Approved for public release, 20-529.

\bibliographystyle{ieeetr}
\bibliography{bibliography}

\begin{thebibliography}{10}

\bibitem{Tompkins}
P.~Tompkins, {\em Mission-directed path planning for planetary rover
  exploration}.
\newblock 2005.

\bibitem{Nieto}
D.~{Nieto-Hern{\'a}ndez}, C.~{M{\'e}ndez-Barrios}, J.~{Escareno},
  V.~{Ram{\'i}rez-Rivera}, L.~A. {Torres}, and H.~{M{\'e}ndez-Az{\'u}a},
  ``Non-holonomic flight modeling and control of a tilt-rotor {MAV},'' in {\em
  2019 6th International Conference on Control, Decision and Information
  Technologies (CoDIT)}, pp.~1947--1952, 2019.

\bibitem{Aguiar}
A.~P. Aguiar and A.~M. Pascoal, ``Regulation of a nonholonomic autonomous
  underwater vehicle with parametric modeling uncertainty using {L}yapunov
  functions,'' in {\em Proceedings of the 40th IEEE Conference on Decision and
  Control}, vol.~5, pp.~4178--4183, IEEE, 2001.

\bibitem{Colombo}
L.~Colombo, R.~Gupta, A.~Bloch, and D.~M. de~Diego, ``Variational
  discretization for optimal control problems of nonholonomic mechanical
  systems,'' in {\em 2015 54th IEEE Conference on Decision and Control (CDC)},
  pp.~4047--4052, IEEE, 2015.

\bibitem{Grushkovaskaya}
V.~Grushkovskaya and A.~Zuyev, ``Obstacle avoidance problem for second degree
  nonholonomic systems,'' in {\em 2018 IEEE Conference on Decision and Control
  (CDC)}, pp.~1500--1505, IEEE, 2018.

\bibitem{Varricchio}
V.~Varricchio and E.~Frazzoli, ``Asymptotically optimal pruning for
  nonholonomic nearest-neighbor search,'' in {\em 2018 IEEE Conference on
  Decision and Control (CDC)}, pp.~4459--4466, IEEE, 2018.

\bibitem{Dubins}
L.~E. Dubins, ``On curves of minimal length with a constraint on average
  curvature, and with prescribed initial and terminal positions and tangents,''
  {\em American Journal of Mathematics}, vol.~79, no.~3, pp.~497--516, 1957.

\bibitem{ReedsShepp}
J.~A. Reeds and L.~A. Shepp, ``Optimal paths for a car that goes both forwards
  and backwards.,'' {\em Pacific J. Math.}, vol.~145, no.~2, pp.~367--393,
  1990.

\bibitem{Barraquand}
J.~Barraquand and J.-C. Latombe, ``Nonholonomic multibody mobile robots:
  Controllability and motion planning in the presence of obstacles,'' {\em
  Algorithmica}, vol.~10, no.~2-4, p.~121, 1993.

\bibitem{AgarwalWang}
P.~K. Agarwal and H.~Wang, ``Approximation algorithms for curvature-constrained
  shortest paths,'' {\em SIAM Journal on Computing}, vol.~30, no.~6,
  pp.~1739--1772, 2001.

\bibitem{Fierro}
R.~Fierro and F.~L. Lewis, ``Control of a nonholonomic mobile robot using
  neural networks,'' {\em IEEE transactions on neural networks}, vol.~9, no.~4,
  pp.~589--600, 1998.

\bibitem{Ali}
Z.~A. Ali, D.~Wang, M.~Safwan, W.~Jiang, and M.~Shafiq, ``Trajectory tracking
  of a nonholonomic wheeled mobile robot using hybrid controller,'' {\em
  International Journal of Modeling and Optimization}, vol.~6, no.~3, p.~136,
  2016.

\bibitem{Deluca}
A.~D. Luca and G.~Oriolo, {\em Modelling and Control of Nonholonomic Mechanical
  Systems}, pp.~277--342.
\newblock Vienna: Springer Vienna, 1995.

\bibitem{Shukla}
A.~Shukla, E.~Singla, P.~Wahi, and B.~Dasgupta, ``A direct variational method
  for planning monotonically optimal paths for redundant manipulators in
  constrained workspaces,'' {\em Robotics and Autonomous Systems}, vol.~61,
  no.~2, pp.~209--220, 2013.

\bibitem{Wu}
W.~Wu, H.~Chen, and P.-Y. Woo, ``Time optimal path planning for a wheeled
  mobile robot,'' {\em Journal of Robotic Systems}, vol.~17, no.~11,
  pp.~585--591, 2000.

\bibitem{Duleba}
I.~Duleba and J.~Z. Sasiadek, ``Nonholonomic motion planning based on newton
  algorithm with energy optimization,'' {\em IEEE transactions on control
  systems technology}, vol.~11, no.~3, pp.~355--363, 2003.

\bibitem{KhanEtAl}
A.~Khan, I.~Noreen, and Z.~Habib, ``An energy efficient coverage path planning
  approach for mobile robots,'' in {\em Intelligent Computing} (K.~Arai,
  S.~Kapoor, and R.~Bhatia, eds.), (Cham), pp.~387--397, Springer International
  Publishing, 2019.

\bibitem{WangEtAl}
T.~Wang, B.~Wang, H.~Wei, Y.~Cao, M.~Wang, and Z.~Shao, ``Staying-alive and
  energy-efficient path planning for mobile robots,'' in {\em 2008 American
  Control Conference}, pp.~868--873, IEEE, 2008.

\bibitem{Verscheure}
D.~{Verscheure}, B.~{Demeulenaere}, J.~{Swevers}, J.~{De Schutter}, and
  M.~{Diehl}, ``Time-energy optimal path tracking for robots: a numerically
  efficient optimization approach,'' in {\em 2008 10th IEEE International
  Workshop on Advanced Motion Control}, pp.~727--732, March 2008.

\bibitem{AlonsoMora}
J.~{Alonso-Mora}, P.~{Beardsley}, and R.~{Siegwart}, ``Cooperative collision
  avoidance for nonholonomic robots,'' {\em IEEE Transactions on Robotics},
  vol.~34, pp.~404--420, April 2018.

\bibitem{Lee2017}
B.~H. Lee, J.~D. Jeon, and J.~H. Oh, ``Velocity obstacle based local collision
  avoidance for a holonomic elliptic robot,'' {\em Autonomous Robots}, vol.~41,
  pp.~1347--1363, Aug 2017.

\bibitem{TakeiTsai1}
R.~Takei, R.~Tsai, H.~Shen, and Y.~Landa, ``A practical path-planning algorithm
  for a simple car: a {H}amilton-{J}acobi approach,'' in {\em Proceedings of
  the 2010 American Control Conference}, pp.~6175--6180, June 2010.

\bibitem{TakeiTsai2}
R.~Takei and R.~Tsai, ``Optimal trajectories of curvature constrained motion in
  the {H}amilton-{J}acobi formulation,'' {\em Journal of Scientific Computing},
  vol.~54, pp.~622--644, Feb 2013.

\bibitem{Parkinson}
C.~Parkinson, D.~Arnold, A.~L. Bertozzi, Y.~T. Chow, and S.~Osher, ``Optimal
  human navigation in steep terrain: a hamilton--jacobi--bellman approach,''
  {\em Communications in Mathematical Sciences}, vol.~17, no.~1, pp.~227--242,
  2019.

\bibitem{Tsitsiklis}
J.~N. Tsitsiklis, ``Efficient algorithms for globally optimal trajectories,''
  {\em IEEE Transactions on Automatic Control}, vol.~40, pp.~1528--1538, Sep
  1995.

\bibitem{Sethian1}
J.~A. Sethian, ``A fast marching level set method for monotonically advancing
  fronts,'' {\em Proceedings of the National Academy of Sciences}, vol.~93,
  no.~4, pp.~1591--1595, 1996.

\bibitem{Bellman}
R.~Bellman, {\em Dynamic Programming}.
\newblock Rand Corporation research study, Princeton University Press, 1957.

\bibitem{CrandallLions}
M.~G. Crandall and P.-L. Lions, ``Viscosity solutions of {H}amilton-{J}acobi
  equations,'' {\em Transactions of the American Mathematical Society},
  vol.~277, no.~1, pp.~1--42, 1983.

\bibitem{Bardi1997}
M.~Bardi and I.~Capuzzo-Dolcetta, {\em Optimal Control and Viscosity Solutions
  of Hamilton-Jacobi-Bellman Equations}.
\newblock Modern Birkh{\"a}user Classics, Birkh{\"a}user Boston, 2008.

\bibitem{TsaiOsherSweep}
Y.-H.~R. Tsai, L.-T. Cheng, S.~Osher, and H.-K. Zhao, ``Fast sweeping
  algorithms for a class of {H}amilton--{J}acobi equations,'' {\em SIAM Journal
  on Numerical Analysis}, vol.~41, no.~2, pp.~673--694, 2003.

\bibitem{CrandallMajda}
M.~G. Crandall and A.~J. Majda, ``Monotone difference approximations for scalar
  conservation laws.,'' 1980.

\bibitem{OsherShu}
S.~Osher and C.-W. Shu, ``High order essentially non--oscillatory schemes for
  {H}amilton--{J}acobi equations,'' {\em SIAM Journal of Numerical Analysis},
  vol.~28, pp.~907--922, August 1991.

\bibitem{Oberman}
A.~M. Oberman, ``Convergent difference schemes for degenerate elliptic and
  parabolic equations: {H}amilton--{J}acobi equations and free boundary
  problems,'' {\em SIAM Journal on Numerical Analysis}, vol.~44, no.~2,
  pp.~879--895, 2006.

\bibitem{Lin}
A.~T. Lin, Y.~T. Chow, and S.~J. Osher, ``A splitting method for overcoming the
  curse of dimensionality in {H}amilton--{J}acobi equations arising from
  nonlinear optimal control and differential games with applications to
  trajectory generation,'' {\em Communications in Mathematical Sciences},
  vol.~16, 1 2018.

\end{thebibliography}





%
%

\end{document}